\documentclass[12pt,reqno,intlim]{amsart}

\usepackage{amssymb,amsmath}

\usepackage[arrow,matrix,curve]{xy}

\newdir{ >}{{}*!/-5pt/\dir{>}}

\newcommand{\nc}{\newcommand}

\nc{\bC}{\bold{C}} \nc{\bN}{\Bbb{N}} \nc{\cF}{\mathcal{F}}
\nc{\cE}{\mathcal{E}} \nc{\cR}{\mathcal{R}} \nc{\cM}{\mathcal{M}}
\nc{\al}{\alpha} \nc{\bt}{\beta} \nc{\gm}{\gamma} \nc{\dl}{\delta}
\nc{\om}{\omega} \nc{\sg}{\sigma} \nc{\Sg}{\Sigma} \nc{\vf}{\varphi}
\nc{\ve}{\varepsilon} \nc{\os}{\overset} \nc{\ol}{\overline}
\nc{\ul}{\underline} \nc{\us}{\underset} \nc{\sbs}{\subset}
\nc{\bsl}{\backslash} \nc{\Ra}{\Rightarrow}
\nc{\lra}{\longrightarrow} \nc{\all}{\allowdisplaybreaks}

\nc{\Codes}{\operatorname{{\bold{Codes}}}}
\nc{\RegMono}{\operatorname{\mathcal{R}{\rm{eg}\mathcal{M}{\rm{ono}\!}}}}
\nc{\RegEpi}{\operatorname{\mathcal{R}{\rm{eg}\mathcal{E}{\rm{pi}\!}}}}
\nc{\Mn}{\operatorname{\mathcal{M}{\rm{ono}\!}}}
\nc{\Ep}{\operatorname{\mathcal{E}{\rm{pi}\!}}}
\nc{\Rg}{\operatorname{\mathcal{R}{\rm{eg}\!}}}
\nc{\Ob}{\operatorname{Ob\!}}

\numberwithin{equation}{section}

\theoremstyle{definition}

\theoremstyle{remark}

\begin{document}

\title[]
{Formulae for the number of non-negative integer solutions of linear Diophantine equation and inequality}

\author{Eteri Samsonadze}
\begin{center}
Proc. Tbilisi State University, 239 (1983), 36-42. 

\end{center}
\begin{center}
Translation from Russian
\end{center}
\begin{abstract} New formulae are presented for the number $P(b)$ of non-negative integer solutions of a Diophantine equation $\sum_{i=1}^{n}a_ix_i=b$ and for the number $Q(b)$ of non-negative integer solutions of the Diophantine inequality $\sum_{i=1}^{n}a_ix_i\leq b$ $(a_i>0, b\geq 0$).\end{abstract}

\maketitle

The methods of finding the number of integer non-negative solutions of a linear Diophantine equation are known (\cite{E}, \cite{A}, \cite{F}, \cite{S}, etc.). In the present paper, employing a new method based on the properties of the Kronecker symbol, yet another formula for the number of integer non-negative solutions of a linear Diophantine equation with an arbitrary number of unknowns is obtained.

Consider an equation
\begin{center}
$$\sum_{i=1}^{n}a_ix_i=b, \ \ \ \ \ \ \ \ \ \ \ \ \ \ \ \ \ \ \ \ \ \ \ \ \ \ \ \ \ \ \ \ \ \ \ \ \ \ \ \ \ \ \ \ \ \ \ \ \ \ \ (1)$$
\end{center}
\noindent where $a_i$ $(i=1,2,...,n)$ are natural numbers; $b$ is an integer non-negative number; $n\geq a$. 

We will prove that, for the number $P(b)$ of integer non-negative solutions of equation (1), the following formula is valid:
\begin{center}
$$P(b)=\sum_{t_1=0}^{d_1-1}\sum_{t_2=0}^{d_2-1}...\sum_{t_n=0}^{d_n-1}C(f(t_1,t_2,...,t_n)+1;n-1),\ \ \ \ \  (2)$$
\end{center}
\noindent where 
$$f(t_1,t_2,...,t_n)=\frac{1}{M}(b-a_1t_1-a_2t_2-...-a_nt_n);$$
$$d_i=\frac{M}{a_i} (i=1,2,...,n);$$
\noindent $M$ is the least common multiple of the numbers $a_1,a_2,...,a_n$;
$$C(k,l)=\frac{1}{l!}k(k+1)...(k+l-1)$$
\noindent if $k$ is a natural number, and $$C(k,l)=0$$
\noindent  otherwise.

\begin{proof}
It is not hard to see that
\begin{center}
$$P(b)=\sum_{x_1}\sum_{x_2}...\sum_{x_n}\delta(\sum_{i=1}^{n}a_ix_i;b),\ \ \ \ \ \ \ \ \ \ \ \ \ \ \ \ \ \ \ \ \ \ \ \ \ \ \ \ \ \  (3)$$
\end{center}
\noindent where $\delta(x;y)$ is the Kronecker symbol, while $x$ in the symbol $\sum_{x}$ takes all non-negative integer values.

Taking formula (3) and the following equality into account
$$\sum_{x}\delta(x;a)=\sum_{x}\sum_{t=0}^{d-1}(dx+t;a),$$
\noindent for any $a$ and natural $d$, we obtain
$$P(b)=\sum_{x_1}\sum_{x_2}...\sum_{x_n}\sum_{t_1=0}^{d_1-1}\sum_{t_2=0}^{d_2-1}...\sum_{t_n=0}^{d_n-1}\delta(\sum_{i=1}^{n}a_i(d_ix_i+t_i);b).$$ 
Since $a_id_i=M$ $(i=1,2,...,n)$, this implies:
$$P(b)=\sum_{x_1}\sum_{x_2}...\sum_{x_n}\sum_{t_1=0}^{d_1-1}\sum_{t_2=0}^{d_2-1}...\sum_{t_n=0}^{d_n-1}\delta(\sum_{i=1}^{n}x_i;\frac{b-\sum_{i=1}^{n}a_it_i}{M}).$$
This formula and formula (3) imply that
\begin{center}
$$P(b)=\sum_{t_1=0}^{d_1-1}\sum_{t_2}^{d_2-1}...\sum_{t_n=0}^{d_n-1}\sum_{t_n=0}^{d_n-1}P_{*}(f(t_1,t_2,...,t_n)),\ \ \ \ \ \ \ \ \ \ \ \ \ \ (4)$$
\end{center}
\noindent where $P_{*}(m)$ denotes the number of integer non-negative solutions of the equation
$$x_1+x_2+...+x_n=m.$$ 

Applying the well-known formula \cite{G}
$$\sum_{x=1}^{k}x(x+1)...(x+l)=\frac{1}{l+2}k(k+1)...(k+l+1),$$
and the principle of the mathematical induction by $n$, it is not hard to prove that
\begin{center}
$$P_{*}(m)=\frac{1}{(n-1)!}(m+1)(m+2)...(m+n-1) \ \ \ \ \ \ \ \ \ \ \ \ \  \ \ \ \ (5)$$
\end{center}
\noindent for $n\geq 2$ and integer $m\geq0$. Formulas (4) and (5) imply (2). 
\end{proof}

\vskip+5mm
Note that the number of summands in the right-hand part of formula (2) does not depend on $b$. At that, if $$b\geq nM-(a_1+a_2+...+a_n),$$ all $f(t_1,t_2,...,t_n)$ are non-negative. Note also that formula (2) is valid also in the case where $M$ is an arbitrary common multiple of the numbers $a_1,a_2,...,a_n$.

Since the number $Q(b)$ of integer non-negative solutions of an inequality 
\begin{center}
$$\sum_{i=1}^{n}a_ix_i\leq b,$$
\end{center}
\noindent is equal to the number of integer non-negative solutions of the equation
$$x_{n+1}+\sum_{i=1}^{n}a_ix_i=b,$$
 formula (2) implies that, for natural $a_i$ $(i=1,2,...,n)$ and integer $b\geq 0$, one has
\begin{center}
$$Q(b)=\sum_{t_1=0}^{d_1-1}\sum_{t_2=0}^{d_2-1}...\sum_{t_n=0}^{d_n-1}\sum_{t_{n+1}=0}^{M-1}C(f(t_1,t_2,...,t_n)+1;n),\ \ \ \ \ \ \ \ \ \ (6)$$
\end{center}
\noindent where $$f(t_1,t_2,...,t_{n+1})=\frac{1}{M}(b-\sum_{i=1}^{n}a_it_i-t_{n+1});$$
$$d_i=\frac{M}{a_i} \ \ \ \ \ \ \  (i=1,2,...,n);$$
\noindent $M$ is the least common multiple of the numbers $a_1,a_2,...,a_n$.
\vskip+3mm
As an implication of formula (2), we give the following formula for calculating $P(b)$:
\begin{center}
$$P(b)=\sum_{i=1}^{n}l_iC(b'+2-i; n-1),\ \ \ \ \ \ \ \ \ \ \ \ \ \ \ \ \ \ \ \ \ \ \ \ \ \  (7)$$
\end{center}
\noindent where 

$b'=\left[\frac{b}{M}\right]$ is the integral part of the number $\frac{b}{M}$;
\vskip+1mm
$l_i=P'(r+(i-1)M)$ $(i=1,2,...,n)$;
\vskip+1mm
$r$ is the remainder of $b$ modulo $M$;
\vskip+1mm
$P'(c)$ is the number of integer non-negative solutions of the system
\vskip+1mm
\vskip+1mm
\vskip+1mm
\begin{center}
$$ \left\{
\begin{array}{llllllll}
\sum_{j=1}^{n}a_jt_j=c,\\
\\
0\leq t_j\leq \frac{M}{a_j}-1 ~ (j=1,2,...,n),
\end{array}
\right.\ \ \ \ \ \ \ \ \ \ \ \ \ \ \ \ \ \ \ \ \ \ \ \ (8)$$
\end{center}
\vskip+3mm
Note that  $l_1$, $l_2$,..., $l_n$ depend only on the numbers $a_1$, $a_2$, ...., $a_n$ and $r$. 
\vskip+3mm
Thus, formula (7) enables us to calculate the number $P(b)$ of integer non-negative solutions of equation (1) for an arbitrary $b$ that is congruent to $r$ modulo $M$ provided that we know the number $P'(c)$ of integer non-negative solutions of system (8) for $c=r$, $c=r+M$, ..., $c=r+(n-1)M$ $(0\leq r <M)$.

Note also that $P'(c)\leq P(c)$ for any $c$ and $P'(c)=0$ for $$c> nM-(a_1+a_2+...+a_n).$$ For finding all $P'(c)$'s it suffices to know the values of the expression $\sum_{j=1}^{n}a_jt_j$ for $c\leq t_j\leq \frac{M}{a_j}-1$ $(j=1,2,...,n)$.
\vskip+3mm
\textbf{Example}. Find the number $P(b)$ of integer non-negative solutions of the equation
$$x_1+x_2+...+x_{n-1}+2x_n=b,$$
\noindent for an integer $b\geq 0$.

To calculate $P(b)$, we will apply formula (7), where we set $$a_1=a_2=...=a_{n-1}=1, $$ $$a_n=2, M=2, t_n=0, 0\leq t_i \leq 1 \ \ (i=1,2,...,n-1).$$ 

It is not hard to see that if $t_1$, $t_2$, ..., $t_{n-1}$ take values $0$ and $1$, then the expression $$t_1+t_2+...+t_{n-1}+2t_n$$ takes the value $0$ precisely one time; the value $1$ -- precisely $C^{1}_{n-1}$ times, etc., the value $n-1$ -- precisely $C_{n-1}^{n-1}$ times. Moreover, $$t_1+t_2+...+t_{n-1}+2t_n\leq n-1.$$ Therefore, $$P'(0)=C_{n-1}^{0}, P'(1)=C^{1}_{n-1},...,P'(n-1)=C_{n-1}^{n-1}$$ and $$P'(c)=0$$ for $c>n-1$.

This fact together with formula (7) implies that
$$P(b)=\sum_{k=0}^{\left[\frac{n-1}{2}\right]} C_{n-1}^{2k}C(b'-k+1;n-1)$$
\noindent for an even $b$ and 
$$P(b)=\sum_{k=0}^{\left[\frac{n-2}{2}\right]} C_{n-1}^{2k+1}C(b'-k+1;n-1)$$
\noindent for an odd $b$, where 
$$C(m;n-1)=\frac{1}{(n-1)!}m(m+1)...(m+n-2)$$
\noindent  if $m>0$, and $$C(m;n-1)=0$$ otherwise; $b'=\left[\frac{b}{2}\right]$.
\vskip+3mm

\vskip+3mm

\vskip+3mm

Received 15.05.1982. 
\vskip+2mm
Tbilisi State University, Department of General Mathematics; 

e-mail: eteri.samsonadze@outlook.com.




\end{document}